\documentclass[12pt]{amsart}

\usepackage[centertags]{amsmath}
\usepackage{amsfonts}
\usepackage{amscd}
\usepackage[sc,osf,slantedGreek]{mathpazo}
\usepackage{euscript}
\usepackage{graphics}
\usepackage{color}
\usepackage[all]{xy}

\newcommand{\B}[1]{{\mathbb #1}}

\newtheorem{theorem}[subsection]{Theorem}%[section]

\newtheorem{corollary}[subsection]{Corollary}

\theoremstyle{definition}

\newtheorem{example}[subsection]{Example}
\theoremstyle{remark}
\newtheorem{remark}[subsection]{Remark}

\numberwithin{figure}{section}
\numberwithin{table}{section}
\newcommand{\Mo}{(M,\omega )}

\newcommand\Ham{\operatorname{Ham}}

\begin{document}

\title{The Lalonde-McDuff conjecture and the fundamental group}
\author{Jarek K\c{e}dra}
\address{Mathematical Sciences, University of Aberdeen, 
         Meston Building, Aberdeen, AB243UE, Scotland, UK\\
Institute of Mathematics, University of Szczecin,
Wielkopolska 15, 70-451 Szczecin, Poland}
%\curraddr{      }
\email{kedra@maths.abdn.ac.uk}
\urladdr{http://www.maths.abdn.ac.uk/\~{}kedra}
\date{\today }
%\thanks{thanks} 
\begin{abstract}
We give a simple proof of the Lalonde-McDuff conjecture for aspherical
manifolds.
\end{abstract}

\maketitle

\section{Introduction}
Let $G$ be a topological group acting on a manifold $M.$ 
Let 
$$
M\stackrel{i}\to M_G \stackrel{\pi}\to BG
$$
be the universal fibration associated to the action.
It is a fundamental question to determine the cohomology
$H^*(M_G)$. This cohomology is also
known as the equivariant cohomology of $M$ associated with
the action of $G$ and denoted by $H_G(M).$ 
A cohomology class in $H^*(M_G)$ is called a $G$-equivariant
class of $M$.

It follows from the Leray-Hirsch theorem (Theorem 4D.1 in Ha\-tcher ~\cite{MR1867354})
that $H(M_G;\B Q)$ is
isomorphic as an $H^*(BG;\B Q)$-module to the tensor product of the
cohomology of the base ans the cohomology of the fibre if and only if
the homomorphism $i^*:H^*(M_G;\B Q)\to H^*(BG;\B Q)$ induced by the
inclusion of the fibre is surjective.

In the present note we make a simple observation about the above
homomorphism under some hypothesis on the group action. As a consequence we
prove the so called Lalonde-McDuff conjecture for aspherical manifolds
(see Section \ref{S:lm} for details).

%%%%%%%%%%%%%%%%%%%%%%%%%%%%%%%%%%%%%%%%%%%%%%%%%%%%%
\section{The main observation}\label{S:main}
%%%%%%%%%%%%%%%%%%%%%%%%%%%%%%%%%%%%%%%%%%%%%%%%%%%%%

Let $c:M\to B\pi_1(M)$ be the map classifying the universal cover.
A cohomology class $\alpha\in H^*(M)$ is called {\em a $\pi_1$-class} if it
is in the image of the induced homomorphism 
$$c^*:H^*(\pi_1(M))\to H^*(M).$$

Let $p\in M$ be a fixed point and let $\operatorname{ev}:G\to M$ be the 
corresponding evaluation map defined by
$$\operatorname{ev}(f) := f(p).$$

\begin{theorem}\label{T:main}
Suppose that the evaluation map induces the trivial homomorphism on
the fundamental group. Then every $\pi_1$-class is the image of some
$G$-equivariant class.
\end{theorem}

\begin{proof}
If the evaluation map induces the trivial homomorphism on the fundamental
group then the connecting homomorphism 
$$
\partial:\pi_2(BG)\to \pi_1(M)
$$ is also trivial. In fact, it is the same homomorphism as the one
induced by the evaluation, after the indentification $\pi_1(G)\cong \pi_2(BG).$
Since $G$ is connected, $BG$ is simply connected and, by the long exact
sequence of homotopy groups of the universal fibration, we get the
isomorphism $i_*:\pi_1(M)\to \pi_1(M_G)$. Thus the classifying map
$c:M\to B\pi_1(M)$ factors through $M_G$ which finishes the proof.
\end{proof}

\begin{corollary}\label{C:aspherical}
Let $M$ be an aspherical manifold. If the evaluation map induces the
trivial homomorphism on the fundamental group then the universal
fibration $M\to M_G\to BG$ is homotopy trivial.
In particular, the $G$-equivariant cohomology of $M$ is isomprphic {\em
as a ring} to the tensor product $H^*(BG)\otimes H^*(M).$
\end{corollary}

\begin{proof}
Since $M$ is aspherical we have that $M=B\pi_1(M).$ 
The map $M_G\to B\pi_1(M)$ classifying the universal cover
is a homotopy inverse of the inclusion
of the fibre $i:M\to M_G.$ This implies that the fibration
is homotopy trivial.
\end{proof}

\begin{example}
It is easy to see that the image of the homomorphism
induced by the evaluation map is contained in the centre of the
fundamental group of $M$. Hence the above results
apply to manifolds whose fundamental group has trivial centre.
In particular, if $M$ is such an aspherical manifold 
then any bundle with fibre $M$
over simply connected base is homotopy trivial.
\end{example}

\begin{remark}
The last observation also follows from the classical result
of Gottlieb \cite[Theorem III.2]{MR32:6454}
which states that the identity component of the
space of homotopy equivalences of an aspherical manifold $M$
is itself aspherical and has the fundamental group isomorphic
to the centre of $\pi_1(M)$.
\end{remark}

\section{The Lalonde-McDuff conjecture}\label{S:lm}
Let $\Mo$ be a symplectic manifold. A locally trivial bundle
$$\Mo\stackrel{i}\to E\stackrel{\pi}\to B$$ is called {\bf
Hamiltonian} if its structure group is a subgroup of the group
$\Ham\Mo$ of Hamiltonian diffeomorphisms of $\Mo.$ La\-lon\-de and
McDuff conjectured in \cite{MR1941438} that the rational cohomology of
the total space is isomorphic as a vector space to the tensor product
of the cohomology of the base and the cohomology of the fibre.  The
conjecture has been proved in many particular cases.  The examples
include fibrations where the fibre is a K\"ahler manifold, fibrations
for which the structure group is a compact Lie group, four dimensional
manifolds \cite{MR1941438} or nilmanifolds \cite{zosia}.  If $\Mo$ is
aspherical we get a stronger statement.

The group of Hamiltonian diffeomorphisms of a closed symplectic manifold
has the property that the evaluation map induces the trivial homomorphism
of the fundamental group. The proof of this fact is nontrivial and can
be found in McDuff-Salamon \cite[Corollary 9.1.2]{MR2045629}. Thus we can
apply Theorem \ref{T:main} in the Hamiltonian case and Corollary \ref{C:aspherical}
gives the following result.

\begin{theorem}\label{T:aspherical}
Let $\Mo\to E\to B$ be a Hamiltonian fibration.
If $M$ is aspherical then the fibration is homotopy trivial.
In particular, the Lalonde-McDuff conjecture holds for 
aspherical manifolds.\qed
\end{theorem}

Combining our argument with some known results we obtain
the following more general theorem.

\begin{theorem}\label{T:lm}
Let $\Mo$ be a symplectic manifold.
If the cohomology ring $H^*(M;\B Q)$ is generated by 
$H^2(M)$ and the Chern classes and the
$\pi_1(M)$-classes then the Lalonde-McDuff conjecture
holds for $\Mo.$
\end{theorem}
\begin{proof}
The fact that the second cohomology is conatined in the
image of the map induced by the inclusion of the fibre
is proven in Lalonde-McDuff \cite[Theroem 1.1]{MR1941438}.

The Chern classes are in the image of $i^*:H^*(E)\to H^*(M)$
because the tangent bundle $TM$ is the pull-back of the
bundle tangent to the fibers of $E\to B.$

Now the statement follows from Theorem \ref{T:main} and the 
Leray-Hirsch theorem.
\end{proof}

%\nocite{*}
\bibliography{../../bib/bibliography}
\bibliographystyle{plain}

\end{document}